\documentclass[12pt]{amsart}

\newtheorem{theorem}{Theorem}[section]

\newtheorem{proposition}[theorem]{Proposition}
\newtheorem{definition}{Definition}[section]

\numberwithin{equation}{section}
\author{Jun Sun}

\address{Jun Sun, School of Mathematical Sciences\\
Wuhan University\\ Wuhan 430072, P. R. of China.}
\email{sunjun@whu.edu.cn}


\thanks {}

\begin{document}

\title[Lagrangian $L$-stability of Lagrangian Translating Solitons]
{Lagrangian $L$-stability of Lagrangian Translating Solitons}

\begin{abstract}
In this paper, we prove that any Lagrangian translating soliton is Lagrangian $L$-stable.
\end{abstract}

\maketitle


\section{Introduction}

\allowdisplaybreaks

\vspace{.1in}

\noindent Recent years, motivated by the problem of existence of special Lagrangian submanifolds, Lagrangian mean curvature flow has attracted much attention. It was proved by Chen-Li (\cite{CL}) and Wang (\cite{Wang}) that there is no finitie time Type I singularity for almost calibrated Lagrangian mean curvature flow. Therefore, there are many works concentrating on Type II singularities of Lagrangian mean curvature flow, especially, on Lagrangian translating solitons (\cite{HL}, \cite{HanS}, \cite{NT}, \cite{Sun1}, \cite{Sun2}, etc.).

\vspace{.1in}

An $n$-dimensional submanifold $\Sigma^n$ in ${\mathbb R}^{n+k}$ is called a {\em translating soliton} if there exists a constant vector ${\bf T}\in {\mathbb R}^{n+k}$, such that
\begin{equation}\label{e-TS}
{\bf T}^{\perp}={\bf H}
\end{equation}
holds on $\Sigma$, where ${\bf H}$ is the mean curvature vector of $\Sigma^n$ in ${\mathbb R}^{n+k}$.

Similar to that of self-shrinkers (\cite{CM}), one can also study the translating solitons from variational viewpoint. Actually, translating solitons can be viewed as critical points of the following functional:
\begin{equation}\label{e-F}
F(\Sigma)=\int_{\Sigma}e^{\langle \bf T,\bf x\rangle}d\mu,
\end{equation}
where $\bf x$ is the position vector in ${\mathbb R}^{n+k}$ and $d\mu$ is the induced area element on $\Sigma$. Then it is natural to define stability of translating solitons. Shahriyari (\cite{Sh}) proved that any translating graph in ${\mathbb R}^3$ is $L$-stable. 

A translating soliton $\Sigma^n$ in ${\mathbb C}^{n}$ is called a {\em Lagrangian translating soliton} if it is also a Lagrangian submanifold of ${\mathbb C}^{n}$. In \cite{Yang}, L. Yang proved that any  Lagrangian translating soliton is Hamiltonian $L$-stable. In this paper, we prove that it is in fact Lagrangian $L$-stable:

\vspace{.1in}

\begin{theorem} \label{thm1.1}
Any Lagrangian translating soliton is Lagrangian $L$-stable.
\end{theorem}

\vspace{.1in}

The proof of Theorem \ref{thm1.1} relies crucially on that the variation is Lagrangian. There are many examples for Lagrangian translating solitons (\cite{CasL}, \cite{JLT}, etc.). By Theorem \ref{thm1.1}, they are all Lagrangian $L$-stable. One natural question is whether we can find examples which are in fact $L$-stable (not just only Lagrangian $L$-stable).
 In \cite{AS}, we showed that the Grim Reaper cylinder $\Gamma\times {\mathbb R}^{n-1}$ is $L$-stable in ${\mathbb R}^{n+1}$, where $\Gamma$ is the Grim Reaper in the plane. This is actually true for any mean convex translating soliton $\Sigma^n$ in ${\mathbb R}^{n+1}$. In this paper, we will show that:

\begin{theorem}\label{thm1.2}
The Lagrangian Grim Reaper cylinder $\Gamma\times {\mathbb R}$ in ${\mathbb C}^2$ is $L$-stable.
\end{theorem}

\vspace{.1in}

For the relations between Lagrangian $F$-stable self-shrinkers and Hamiltonian $F$-stable self-shrinkers, we refer to \cite{LZ}. 

\vspace{.1in}

\textbf{Acknowledgement:} The author was supported by the National Natural Science Foundation of China, No. 11401440. Part of the work was finished when the author was a visiting scholar at Massachusetts Institute of Technology (MIT) supported by China Scholarship Council (CSC). The author would like to express his gratitude to Professor Tobias Colding for his invitation, to MIT for their hospitality, and to CSC for their support.

\vspace{.2in}

\section{Preliminaries}

\vspace{.1in}

\noindent In this section, we will recall some results for the first variation and second variation formulas. Since the proofs can be found in Section 4 of \cite{AS} with $f=\langle \bf T,\bf x\rangle$, where we dealt with more general cases (see also \cite{Yang}), we omit the details here.

\vspace{.1in}

Recall that the $F$-functional is defined by 
\begin{equation*}
F(\Sigma)=\int_{\Sigma}e^{\langle \bf T,\bf x\rangle}d\mu.
\end{equation*}
The following first variation formula is known (Propositon 4.1 of \cite{AS}):

\vspace{.1in}

\begin{proposition}\label{prop2.1}
Let $\Sigma_s^n\subset {\mathbb R}^{n+k}$ be a smooth compactly supported variation of $\Sigma$ with normal variational vector field ${\bf V}$, then
 \begin{equation}\label{e2.1}
\frac{d}{ds}|_{s=0}F(\Sigma_s)=\int_{\Sigma}\langle {\bf T}^{\perp}-{\bf H},{\bf V}\rangle e^{\langle \bf T,\bf x\rangle}d\mu.
\end{equation}
In particular, $\Sigma$ is a critical point of $F$ if and only if ${\bf T}^{\perp}={\bf H}$, i.e., $\Sigma$ is a translating soliton in ${\mathbb R}^{n+k}$.
\end{proposition}

\vspace{.1in}

For the second variation formula, we have (see (4.17) of \cite{AS}):

\begin{theorem}\label{thm2.2}
Suppose that $\Sigma$ is a critical point of $F$. If $\Sigma_s^n\subset {\mathbb R}^{n+k}$ be a smooth compactly supported variation of $\Sigma$ with normal variational vector field ${\bf V}$,
then the second variation formula is given by
 \begin{equation}\label{e2.2}
F'':=\frac{d^2}{ds^2}|_{s=0}F(\Sigma_s)=-\int_{\Sigma}\langle L{\bf V},{\bf V}\rangle e^{\langle \bf T,\bf x\rangle}d\mu.
\end{equation}
Here, the stability operator $L$ is defined on a normal vector field $\textbf{V}$ on $M$ by
\begin{equation}\label{e2.3}
  L\textbf{V}=\left(\Delta V^{\alpha}+\langle \textbf{T}, \nabla V^{\alpha}\rangle
              +g^{ik}g^{jl}h^{\alpha}_{ij}h^{\beta}_{kl}V^{\beta}\right) e_{\alpha},
\end{equation}
where $\{e_{\alpha}\}$ is a local orthonormal frame of the normal bundle $N\Sigma$, $g_{ij}$ is the induced metric on $\Sigma$ and ${\bf V}=V^{\alpha}e_{\alpha}$.
\end{theorem}

\vspace{.1in}

\begin{definition}
A translating soliton $\Sigma^n$ in ${\mathbb R}^{n+k}$ is said to be {\bf $L$-stable} if for every compactly supported normal variational vector field ${\bf V}$, we have 
$$F''=-\int_{\Sigma}\langle L{\bf V},{\bf V}\rangle e^{\langle \bf T,\bf x\rangle}d\mu\geq 0.$$
\end{definition}

\vspace{.1in}

Now we turn to Lagrangian translating solitons. Let $\bar \omega$ and $J$ be the standard K\"ahler form and complex structure on ${\mathbb C}^n$, respectively. A  submanifold $\Sigma^n$ is said to be a {\em Lagrangian} submanifold of ${\mathbb C}^n$, if $\bar\omega|_{\Sigma}=0$, or eqivalently, $J$ maps the tangent space of $\Sigma$ on to its normal space at each point of $\Sigma$. For a Lagrangian submanfiold, there is a canonical correspondence between the sections of the normal bundle and the space of 1-forms on $\Sigma$:
\begin{eqnarray*}
\Gamma(N\Sigma) & \longrightarrow & \Lambda^1(\Sigma)\\
\bf V &\longleftrightarrow & \theta_{\bf V}:=-i_{\bf V}\bar\omega.
\end{eqnarray*}
A normal vector field $\bf V$ is a {\em Lagrangian variation} if $\theta_{\bf V}$ is closed; a normal vector field $\bf V$ is a {\em Hamiltonian variation} if $\theta_{\bf V}$ is exact.

\vspace{.1in}

\begin{definition}
A Lagrangian translating soliton $\Sigma^n$ in ${\mathbb C}^{n}$ is said to be {\bf Lagrangian (resp. Hamitlonian) $L$-stable} if for every compactly supported normal Lagrangian (resp. Hamiltonian) variation ${\bf V}$, we have 
$$F''=-\int_{\Sigma}\langle L{\bf V},{\bf V}\rangle e^{\langle \bf T,\bf x\rangle}d\mu\geq 0.$$
\end{definition}

\vspace{.2in}

\section{Lagrangian $L$-stability of Lagrangian Translating Solitons}

\vspace{.1in}

\noindent In this section, we will prove Theorem 1.1. First, we would like to rewrite the second variation formula for Lagrangian variations.

\vspace{.1in}

Let $({\mathbb C}^n,\bar g,J,\bar\omega)$ be the complex Euclidean space with standard metric $\bar g$, complecx structure $J$ and K\"ahler form $\bar\omega$  such that $\bar g=\bar\omega(\cdot,J\cdot)$. Given any Lagrangian submanifold $\Sigma^n$ in ${\mathbb C}^n$, we choose a local orthonormal frame $\{e_i\}_{i=1}^n$ of $T\Sigma$, and set $\nu_i=Je_i$. Then $\{\nu_i\}_{i=1}^n$ forms a local orthonormal frame of the normal bundle $N\Sigma$. The frame can be chosen so that at a fixed point $x\in\Sigma$, we have $\nabla_{e_i}e_j=0$, where $\nabla$ is the induced connection on $\Sigma$. The second fundamental form is defined by
\begin{equation*}
h_{ijk}=\bar g(\overline\nabla_{e_i}e_j,\nu_k),
\end{equation*}
which is symmetric in $i$, $j$ and $k$. Th mean curvature vector is given by
\begin{equation*}
{\bf H}=H_k\nu_k=h_{iik}\nu_k.
\end{equation*}

Let $\{\omega^i\}_{i=1}^{n}$ be the dual basis of $\{e_i\}_{i=1}^{n}$. Then for any normal vector field ${\bf V}=V_i\nu_i$, we have the correspondence 
\begin{equation*}
\theta_{\bf V}:=-i_{\bf V}\bar\omega=V_i\omega^i.
\end{equation*}
Since $d\theta_{\bf V}=\nabla_{e_i}V_j \omega^j\wedge \omega^i$, we see that

\begin{proposition}\label{prop3.1}
A normal vector field $\bf V$ of a Lagrangian submanfiold $\Sigma^n$ in ${\mathbb C}^n$ is a Lagrangian variation if and only if $\nabla_{e_i}V_j=\nabla_{e_j}V_i$. 
\end{proposition}

\vspace{.1in}

Using the above notations, we see that the stability operator (\ref{e2.3}) can be rewritten as
\begin{equation}\label{e3.1}
  L\textbf{V}=\left(\Delta V_i+\langle \textbf{T}, \nabla V_i\rangle
              +h_{kli}h_{klj}V_j\right) \nu_i.
\end{equation}
Therefore, we have

\begin{proposition}\label{prop3.2}
A Lagrangian translating soliton $\Sigma^n$ in ${\mathbb C}^{n}$ is Lagrangian $L$-stable if and only if for every compactly supported normal Lagrangian variation ${\bf V}=V_i\nu_i$, we have 
\begin{equation}\label{e3.2}
F''=-\int_{\Sigma}\left(V_i\Delta V_i+V_i\langle \textbf{T}, \nabla V_i\rangle
              +h_{kli}h_{klj}V_iV_j\right) e^{\langle \bf T,\bf x\rangle}d\mu\geq 0.
\end{equation}
\end{proposition}

Now, we can prove the first main result. For the purpose of convenience, we rewrite it here:

\vspace{.1in}

\begin{theorem} \label{thm3.3}
Any Lagrangian translating soliton is Lagrangian $L$-stable.
\end{theorem}

\vspace{.1in}

{\textbf Proof:} By Proposition \ref{prop3.2}, it suffices to prove that (\ref{e3.2}) holds for every compactly supported Lagrangian variation ${\bf V}=V_i\nu_i$. Since ${\bf V}=V_i\nu_i$ is Lagrangian, by Proposition \ref{prop3.1}, we see that $\nabla_{e_i}V_j=\nabla_{e_j}V_i$. By Ricci identity, we have
\begin{equation}\label{e3.3}
\Delta V_i=\nabla_j\nabla_jV_i=\nabla_j\nabla_iV_j=\nabla_i\nabla_jV_j+R_{jijk}V_k=\nabla_i\nabla_jV_j+R_{ik}V_k,
\end{equation}
where $R_{ik}$ is the Ricci curvature of the induced metric on $\Sigma$. By Gauss equation, we have that
\begin{equation*}
R_{ijkl}=h_{pik}h_{pjl}-h_{pil}h_{pjk},
\end{equation*}
which implies that
\begin{equation}\label{e3.4}
R_{ik}=g^{jl}R_{ijkl}=H_ph_{pik}-h_{pji}h_{pjk}.
\end{equation}
Putting (\ref{e3.4}) into (\ref{e3.3}) yields
\begin{equation*}
\Delta V_i=\nabla_i\nabla_jV_j+R_{ik}V_k=\nabla_i\nabla_jV_j+H_ph_{pik}V_k-h_{pji}h_{pjk}V_k.
\end{equation*}
Therefore, we have
\begin{equation}\label{e3.5}
F''=-\int_{\Sigma}\left(V_i\nabla_i\nabla_jV_j+V_i\langle \textbf{T}, e_j\rangle\nabla_j V_i
              +H_ph_{pij}V_iV_j\right) e^{\langle \bf T,\bf x\rangle}d\mu.
\end{equation}
Integrating by part, we can compute the first term on the right hand side of (\ref{e3.5}) as:
\begin{eqnarray}\label{e3.6}
& & -\int_{\Sigma}V_i\nabla_i\nabla_jV_je^{\langle \bf T,\bf x\rangle}d\mu\nonumber\\
&=& \int_{\Sigma}(\nabla_iV_i\nabla_jV_j+V_i\nabla_jV_j \nabla_i\langle \bf T,\bf x\rangle)e^{\langle \bf T,\bf x\rangle}d\mu\nonumber\\
&=&  \int_{\Sigma}\left[\left(\sum_{j=1}^{n}\nabla_jV_j\right)^2+\left(\sum_{j=1}^{n}\nabla_jV_j\right)\left(\sum_{i=1}^{n} \langle {\bf T},e_i\rangle V_i\right)\right]e^{\langle \bf T,\bf x\rangle}d\mu.
\end{eqnarray}
On the other hand, from the translating soliton equation (\ref{e-TS}), we can easily see that $H_p=\langle\bf T,\nu_p\rangle$. Therefore, using the fact that $\nabla_{e_i}V_j=\nabla_{e_j}V_i$, the second term on the right hand side of (\ref{e3.5}) can be computed as:
\begin{eqnarray}\label{e3.7}
& & -\int_{\Sigma}V_i\langle \textbf{T}, e_j\rangle\nabla_j V_ie^{\langle \bf T,\bf x\rangle}d\mu
      =-\int_{\Sigma}(\nabla_i V_j)V_i\langle \textbf{T}, e_j\rangle e^{\langle \bf T,\bf x\rangle}d\mu\\
&=& \int_{\Sigma}\left[(\nabla_iV_i)\langle \textbf{T}, e_j\rangle V_j+V_iV_j \nabla_i\langle {\bf T},e_j\rangle+V_iV_j \langle {\bf T},e_j\rangle \nabla_i\langle \bf T,\bf x\rangle\right]e^{\langle \bf T,\bf x\rangle}d\mu\nonumber\\
&=& \int_{\Sigma}\left[(\nabla_iV_i)\langle \textbf{T}, e_j\rangle V_j+V_iV_j \langle {\bf T},h_{pij}\nu_p\rangle+V_iV_j \langle {\bf T},e_j\rangle\langle {\bf T},e_i\rangle\right]e^{\langle \bf T,\bf x\rangle}d\mu\nonumber\\
&=&  \int_{\Sigma}\left[\left(\sum_{j=1}^{n}\nabla_jV_j\right)\left(\sum_{i=1}^{n} \langle {\bf T},e_i\rangle V_i\right)+H_ph_{pij}V_iV_j+\left(\sum_{i=1}^{n}\langle {\bf T},e_i\rangle V_i\right)^2\right]e^{\langle \bf T,\bf x\rangle}d\mu.\nonumber
\end{eqnarray}
Here, we used the fact that $\overline\nabla_{e_i}e_j=h_{pij}\nu_j$ at a fixed point by the choice of the frame.
Putting (\ref{e3.6}) and (\ref{e3.7}) into (\ref{e3.5}) yields
\begin{equation*}
F''=\int_{\Sigma}\left(\sum_{j=1}^{n}\nabla_jV_j+\sum_{i=1}^{n}\langle {\bf T},e_i\rangle V_i\right)^2 e^{\langle \bf T,\bf x\rangle}d\mu\geq0.
\end{equation*}
This finishes the proof of the theorem.
\hfill Q.E.D.

\vspace{.1in}

\section{The Lagrangian Grim Reaper Cylinder}

\vspace{.1in}

\noindent In the previous section, we proved that any Lagrangian translating soliton is Lagrangian $L$-stable. However, it is not clear that whether they are $L$-stable. In this section, as an example, we will show that the Grim Reaper cylinder $\Gamma\times \mathbb R$ is in fact $L$-stable in ${\mathbb C}^2$.

\vspace{.1in}

First recall that the Grim Reaper $\Gamma$ in the plane is defined by
\begin{eqnarray*}
\gamma: \left(-\frac{\pi}{2},\frac{\pi}{2}\right) & \longrightarrow & \mathbb C\\
                                                               x     & \longrightarrow & \gamma(x)=(-\log \cos x,x).
\end{eqnarray*}
Then the Grim Reaper cylinder $\Gamma\times \mathbb R$ is defined by
\begin{eqnarray*}
\Phi: \left(-\frac{\pi}{2},\frac{\pi}{2}\right)\times{\mathbb R} & \longrightarrow & {\mathbb C}^2\\
                                                               (x,y)     & \longrightarrow & \Phi(x,y)=(-\log \cos x,x,y,0).
\end{eqnarray*}
We will see that it is a Lagrangian translating soliton and is $L$-stable.

\vspace{.1in}

\begin{theorem}\label{thm4.1}
The Grim Reaper cylinder $\Sigma=\Gamma\times \mathbb R$ is a Lagrangian translating soliton of ${\mathbb C}^2$ and is $L$-stable.
\end{theorem}

\vspace{.1in}

{\textbf Proof:} By the definition of $\Phi$, the tangent space of $\Sigma$ is spanned by 
\begin{equation*}
\Phi_x=(\tan x,1,0,0),  \ \ \ \ \Phi_y=(0,0,1,0).
\end{equation*}
The orthonormal basis of the normal space can be taken as
\begin{equation*}
e_3=(\cos x,-\sin x,0,0),  \ \ \ \ e_4=(0,0,0,-1).
\end{equation*}
The induced metric can be represented as
\begin{equation}\label{e4.1}
(g_{ij})_{1\leq i,j, \leq 2}= \left(
  \begin{array}{cc}
    \frac{1}{\cos^2 x} & 0 \\
        0       & 1 \\
  \end{array}
\right),
\ \ \ \
(g^{ij})_{1\leq i,j, \leq 2}= \left(
  \begin{array}{cc}
   \cos^2 x & 0 \\
        0       & 1 \\
  \end{array}
\right).
\end{equation}
The induced area form is given by
\begin{equation}\label{e4.2}
d\mu=\sqrt{\det(g_{ij})} dxdy=\frac{1}{\cos x}dxdy. 
\end{equation}
Since
\begin{equation*}
\Phi_{xx}=(\frac{1}{\cos^2 x},0,0,0),  \ \ \Phi_{xy}=(0,0,0,0), \ \ \Phi_{yy}=(0,0,0,0),
\end{equation*}
from $h^{\alpha}_{ij}=\langle \Phi_{ij},e_{\alpha}\rangle$,  we see that the second fundamental form are given by
\begin{equation}\label{e4.3}
h^3_{xx}=\frac{1}{\cos x},  \ h^3_{xy}=h^3_{yy}=h^4_{xx}=h^4_{xy}=h^4_{yy}=0.
\end{equation}
Therefore, 
\begin{equation*}
H^3=g^{ij}h^3_{ij}=g^{xx}h^3_{xx}=\cos x, \ H^4=g^{ij}h^4_{ij}=0,
\end{equation*}
and the mean curvature vector is given by
\begin{equation*}
{\bf H}=H^3e_3+H^4e_4=\cos x e_3.
\end{equation*}
Now if we take ${\bf T}=(1,0,0,0)\in {\mathbb C}^2$, then
\begin{equation*}
{\bf T}^{\perp}=\langle {\bf T},e_3\rangle e_3+\langle {\bf T},e_4\rangle e_4=\cos x e_3={\bf H}.
\end{equation*}
Therefore, $\Sigma$ is a translating soliton in ${\mathbb C}^2$.

\vspace{.1in}

Recall that the standard complex structure in ${\mathbb C}^2$ is given by
\begin{equation*}
J=\left(
  \begin{array}{cccc}
    0 & 1  &  0 & 0\\
   -1 &  0 &  0 &  0 \\
     0 &  0  & 0 & 1 \\
    0  &   0 & -1 & 0 \\
  \end{array}
\right).
\end{equation*}
Since 
\begin{equation*}
J\Phi_x=(1,-\tan x,0,0)=\frac{1}{\cos x}e_3, \ J\Phi_y=(0,0,0,-1)=e_4,
\end{equation*}
we see that $\Sigma$ is a Lagrangian translating soliton in ${\mathbb C}^2$.

\vspace{.1in}

Next we will show that $\Sigma$ is $L$-stable. Since
\begin{equation*}
\Delta v+\langle {\bf T},\nabla v\rangle=e^{-\langle{\bf T},{\bf x}\rangle}div_{\Sigma}(e^{\langle{\bf T},{\bf x}\rangle}\nabla v),
\end{equation*}
for any smooth function $v$ on $\Sigma$, we can see easily from Theorem \ref{thm2.2} that $\Sigma$ is $L$-stabe if and only if
\begin{equation}\label{e4.4}
 \int_{\Sigma}g^{ik}g^{jl}h^{\alpha}_{ij}h^{\beta}_{kl}V^{\alpha}V^{\beta} e^{\langle \bf T,\bf x\rangle}d\mu 
  \leq  \int_{\Sigma}\sum_{\alpha}|\nabla V^{\alpha}|^2 e^{\langle \bf T,\bf x\rangle}d\mu
\end{equation}
holds for every compactly supported normal variation vector field ${\bf V}=V^{\alpha}e_{\alpha}$.

\vspace{.1in}

In our case, $\langle {\bf T},{\bf x}\rangle=\langle (1,0,0,0),(-\log \cos x,x,y,0)\rangle=-\log \cos x$ so that
\begin{equation}\label{e4.5}
e^{\langle \bf T,\bf x\rangle}=\frac{1}{\cos x}.
\end{equation}
By (\ref{e4.1}) and (\ref{e4.3}), we have
\begin{equation*}
g^{ik}g^{jl}h^{\alpha}_{ij}h^{\beta}_{kl}V^{\alpha}V^{\beta}=g^{xx}g^{xx}h^3_{xx}h^3_{xx}V^3V^3=\cos^2x (V^3)^2.
\end{equation*}
Combining with (\ref{e4.2}) and (\ref{e4.5}), we get that
\begin{equation}\label{e4.6}
 \int_{\Sigma}g^{ik}g^{jl}h^{\alpha}_{ij}h^{\beta}_{kl}V^{\alpha}V^{\beta} e^{\langle \bf T,\bf x\rangle}d\mu 
 =  \int_{-\infty}^{\infty}\int_{-\frac{\pi}{2}}^{\frac{\pi}{2}}\left(V^3\right)^2dxdy.
\end{equation}
On the other hand, using (\ref{e4.1}), we compute
\begin{equation*}
|\nabla V^{\alpha}|^2=g^{ij}\frac{\partial}{\partial x^i}V^{\alpha}\frac{\partial}{\partial x^j}V^{\alpha}
=\cos^2 x \left(\frac{\partial}{\partial x}V^{\alpha}\right)^2+\left(\frac{\partial}{\partial y}V^{\alpha}\right)^2.
\end{equation*}
Therefore,
\begin{eqnarray}\label{e4.7}
 \int_{\Sigma}\sum_{\alpha}|\nabla V^{\alpha}|^2 e^{\langle \bf T,\bf x\rangle}d\mu
 &=& \int_{-\infty}^{\infty}\int_{-\frac{\pi}{2}}^{\frac{\pi}{2}}\left[\left(\frac{\partial}{\partial x}V^{3}\right)^2+\left(\frac{\partial}{\partial x}V^{4}\right)^2\right]dxdy\nonumber\\
 && +\int_{-\infty}^{\infty}\int_{-\frac{\pi}{2}}^{\frac{\pi}{2}}\frac{1}{\cos^2x}\left[\left(\frac{\partial}{\partial y}V^{3}\right)^2+\left(\frac{\partial}{\partial y}V^{4}\right)^2\right]dxdy
\end{eqnarray}
Note that $V^3,V^4\in C_0^{\infty}\left(\left(-\frac{\pi}{2},\frac{\pi}{2}\right)\times{\mathbb R}\right)$. In particular, for each fixed $y$, we have $V^3(\cdot,y)\in C_0^{\infty}\left(\left(-\frac{\pi}{2},\frac{\pi}{2}\right)\right)$. By Wirtinger inequality, we have for each $y\in {\mathbb R}$ that
\begin{equation*}
\int_{-\frac{\pi}{2}}^{\frac{\pi}{2}}\left(V^3(x,y)\right)^2dx\leq \int_{-\frac{\pi}{2}}^{\frac{\pi}{2}}\left(\frac{\partial}{\partial x}V^{3}(x,y)\right)^2dx
\end{equation*}
Integrating with respect to $y$ yields
\begin{equation}\label{e4.8}
\int_{-\infty}^{\infty}\int_{-\frac{\pi}{2}}^{\frac{\pi}{2}}\left(V^3(x,y)\right)^2dxdy\leq \int_{-\infty}^{\infty}\int_{-\frac{\pi}{2}}^{\frac{\pi}{2}}\left(\frac{\partial}{\partial x}V^{3}(x,y)\right)^2dxdy.
\end{equation}
Combining (\ref{e4.6}), (\ref{e4.7}) and (\ref{e4.8}), we see that (\ref{e4.4}) holds for every compactly supported normal variation vector field ${\bf V}=V^3e_3+V^4e_4$. This shows that the Lagrangian translating soliton $\Sigma$ is $L$-stable.
\hfill Q.E.D.

\end{document}